\documentclass[a4paper,11pt]{article}
\usepackage{amsfonts}

\author{Dmytro Taranovsky}
\title{Constructive Mathematical Truth}
\date{June 1, 2006}

\begin{document}
\maketitle
\begin{abstract}
We define constructive truth for arithmetic and for intuitionistic analysis, and investigate its properties.  We also prove that the set of constructively true (first order) arithmetical statements is $\Pi^1_2$ and $\Sigma^1_2$ hard, and we conjecture it to be complete for second order arithmetic.   A statement is constructively true iff it is realized by a constructive function under continuous function realizability.
\end{abstract} 
\section{Introduction}

\subsection{Historical Background}
Being able to construct a mathematical object is more satisfying than simply knowing it exists, and a branch of mathematics was developed with emphasis on construction. However, a satisfactory description of constructive truth has previously eluded mathematicians. Some like L.E.J. Brouwer thought that mathematics is a free activity of the mind, and hence subjective. Some others confined themselves to recursive functions and analysis. Still others emphasized constructive proof (as in Brouwer-Heyting-Kolmogorov semantics) over truth.

Constructive mathematics is as old as mathematics itself, and most mathematics before the 19th century was constructive. Constructive mathematics was considered just ordinary mathematics and was done in classical logic. In the early 20th century, Brouwer pioneered constructive mathematics, specifically intuitionism, as a distinct branch of mathematics with the rules of logic weakened to ensure that non-constructive results are not provable. Specifically, the law of the excluded middle $A \vee \neg A$ was omitted in favor of weaker statements. However, Brouwer adopted a polemic style and insisted on subjectivity of mathematics. As a consequence, his constructive analysis was imprecise and difficult to use. Only later was the intuitionistic logic formalized.  Around 1950s and later, Russian mathematicians developed a different branch of constructive mathematics, namely constructive recursive analysis. In constructive recursive analysis, every real number is given by a recursive code and every answer must be a construction. 

A major development was Bishop's \textit{Foundations of Constructive Analysis} [Bishop 1967], where he showed that much of mathematics can be done constructively and without using controversial or classically incorrect principles. Bishop accomplished that by using mathematically relevant formalizations (among many classically but not constructively equivalent ones) and by using additional assumptions (like uniform continuity) for his theorems.

Still, semantics of constructive mathematics remained unsatisfactory. There is a qualitative leap from simply working in a constructive formal system to knowing what is constructive truth. Semantics gives meaning to symbol manipulations, explains why some rules are accepted while others are not, makes it much easier to have intuition about the subject, and provides a way resolve the incompleteness in the formal system. Historically, because of usability problems and unsatisfactory semantics, relatively few mathematicians actually worked in constructive formal systems. Some semantics was provided by Kripke models and related developments, but these do not provide a notion of constructive truth. Different statements are true at different nodes of a typical Kripke model.

For constructive recursive mathematics, semantics was provided through Kleene's recursive realizability (which is defined in this paper). However, recursive realizability essentially forbids non-recursive functions, which makes it too restrictive for ordinary mathematics. Later [Kleene 1965], Kleene defined function realizability, which comes close to constructive truth. Unfortunately, Kleene's function realizability became relatively obscure, and I did not know its definition until this paper was essentially complete.

\subsection{Motivation and Outline}

My search for constructive truth started from the question, \textit{What do we know about a statement (beyond its truth) from having a constructive proof regardless of the constructive formal system used?} The answer is that we know a constructive witness to the statement. A statement is constructively true iff it has a constructive witness.

Some constructivists may reject the notion that there is constructive arithmetical truth which can be formally defined. For them, our quest can be viewed in this way: If we have to formally define constructive truth, what would be the definition?  Besides, the mathematical results can be enjoyed whether or not the reader accepts our philosophical claim.  

Clearly, we cannot identify constructive truth with constructive proof since there is no single privileged formal notion of constructive proof. Our definition diverges from the thought of Brouwer and others, but only to the extent necessary for basic properties of truth. Constructive truth resembles recursive realizability, except that every constructively true arithmetical statement is true. In fact, our definition of witnessability (which we also call continuous function realizability) agrees with Kleene's function realizability but for notational differences (we use sets rather than functions, and our treatment of witnesses is more explicit).

In retrospect, defining the exact constructive content of implication is the technical impediment to the definition of constructive truth. We have to use real numbers to overcome the difficulty, which leads to the amazingly high complexity of constructive arithmetical truth. However, the technical difficulty was essentially overcome by Kleene, and the problem was more of a philosophical one. Intuitionism as a philosophy rejects objective mathematical truth, so constructive truth was not pursued. The result was a much less usable intuitionistic analysis than it could have been.

The paper is mostly self-contained. Background information can be found in articles ``Intuitionistic Logic'' and ``Constructive Mathematics'' in the \textit{Stanford Encyclopedia of Philosophy} (online); in the chapter ``Aspects of Constructive Mathematics'' in the \textit{Handbook of Mathematical Logic}; and also in [Moschovakis 2005].  The latter defines function realizability.

The second section defines constructive arithmetical truth, and the third section describes its basic properties. In the fourth section, we prove that constructive arithmetical truth is $\Pi^1_2$ and $\Sigma^1_2$ hard.  For ease of understanding, we build up to that result by presenting easier results first.  We also prove that for a large class of statements, constructive truth is $\Pi^1_1$.  In fifth section, we consider a modification of constructive truth and prove that it is $\Pi^1_1$ complete.  In the sixth section, which is independent of the fourth and fifth sections, we define and analyze constructive truth for intuitionistic analysis.  For brevity and smoothness of flow, in that section we omit formal proofs, but given enough details to make the proofs routine.  That section also discusses constructive reverse mathematics and concludes with an independence result.  The final section discusses correctness of the definition of constructive truth.

\section{Definition of Constructive Truth}
A statement is constructively true if it is true and there is a constructive witness to existential aspects of its truthfulness. Logical connectives are used to designate which aspects are existential. Logically equivalent (in the classical sense) statements may have different designations, and thus different requirements on the witness, which can make the statements constructively inequivalent.

A witness must produce a certain output when given a certain input as described below.  A witness is constructive if there is a constructive mechanism for its operation; the standard and default meaning of this is that the witness is computable, but other notions are possible. The stricter the notion of being constructive, the more difficult it is for a statement to be constructively true; the weakest possible notion is for (classical) truth.

For readers who are familiar with Kleene's function realizability, we note that our definition of witnesses agrees up to recursive interconversion with function realizability; and a statement is constructively true iff it is realized by a recursive function.  The agreement extends to constructive analysis.  However, here we present witnesses in a more explicit and readable syntactic form. 

The language consists of the arithmetical symbols `0', `1', `+', `*', `$<$', `=', logical symbols of `$\neg$', `$\wedge$', `$\vee$', `$\rightarrow$', `$\exists$', `$\forall$', and the modal operator `$\Box$' of constructive truth in the recursive sense. A formula is arithmetical if it does not use `$\Box$'.   Quantifiers and `$\Box$' have the least possible scope. Implication has lower precedence than other logical operators.

\begin{itemize}
    \item A witness for $A$ exists iff $A$ is true.
    \item A witness for $\Box A$ must provide a standard code for a recursive witness for $A$.
    \item A witness for $\neg A$ or for an atomic formula is not required to do anything (but only exists if the formula is true).
    \item A witness for $A\vee B$ must choose $A$ or $B$ and provide a witness for it.
    \item A witness for $A\wedge B$ must provide witnesses for both $A$ and $B$.
    \item A witness for $\exists x A(x)$ must provide $n$ and a witness for $A(n)$.
    \item A witness for $\forall x A(x)$ must when given $n$ as input, provide a witness for $A(n)$.
    \item A witness for $A \rightarrow B$ must provide a witness for $B$ whenever given a witness for $A$.
\end{itemize}

Here $A$ and $B$ are formulas with an assignment of free variables.  Technically, witnesses are real numbers (or arbitrary sets of integers) that specify which output to give for each possible input and satisfy the above requirements. The real number codes a sequence of input-output pairs in arbitrary order (and with arbitrary delay or whitespace), possibly omitting incorrect inputs. Having a witness refers to the ability to query the digits of the real number.  $W$ witnesses $A \rightarrow B$ iff for every $X$ witnessing $A$, $W_X$ witnesses $B$. Here $W$ is viewed as a continuous function, and $W_X$ is $W$ applied to $X$ (using whitespace in $W_X$ to make $f$: $f(W, X)=W_X$ total (and recursive); for inappropriate $X$, $W_X$ might consist of just whitespace).  The rest is a formalization of the syntactic details (that should not affect constructive truth).

Input and output are given in the order of decreasing scope; `$\wedge$' and `$\vee$' are treated like quantifiers.  For example, for $\forall x\exists y (x=2*y \vee  x=2*y+1)$, an input-output pair can be $(25: 12, 1)$.  $25$ is the input $x$; the output provides $12$ as $y$, and uses $1$ to choose the second disjunct ($0$ would indicate the first disjunct; analogously with conjuncts).  $(25: 12)$ is also an input-output pair, but with insufficient output.  $(:)$ is the trivial pair (and is correct here); $(: 12)$ gives too much output.  A sample witness for $\forall x \exists y \: y=2*x$ is ``$(:) \: (0:0)\: (1:2)\: (2:4) \ldots$''. Input does not affect output outside the scope.  For example, for $\forall x\exists y\forall z P(x, y, z)$, $y$ is a part of output rather than input and may not depend on $z$.

For $A \rightarrow B$, input consists of an initial segment of a witness for $A$ and input for $B$.  Arbitrary output (of the right syntactic form) for the consequent is permitted for false initial segments (including whenever $A$ is false).  Segments that are too short (or incorrect) can be rejected in that input-output pairs are omitted.  However, for every witness for $C$, for every tuple of witnesses used as input into $C$, for every appropriate numerical input for $C$, some tuple of finite initial segments of the witnesses combined with the numerical input is accepted.

Output is uniquely determined by the input and is preserved if we extend the initial segments of the candidate witnesses for the antecedents (even for false candidate witnesses, where the output can otherwise be arbitrary; also, whitespace matters here for what constitutes an extension).  (Alternatively, we could have allowed multiple outputs as long as the witness follows through and is accountable for each output.  Any such witness can be converted into a witness with unique outputs.)

The semantic content of an input-output pair for a statement is the assertion (in classical logic) that the pair is correct.  For example, for $\forall x \exists y \: y=x+1 \rightarrow \forall x \exists y \: y=x+2$, the input-output pair (``(2:3) (3:4)'', 2: 4) has semantic content $\forall x \exists y \: y=x+1 \wedge 3=2+1 \wedge 4=3+1 \rightarrow 4=2+2$.

Equivalently, witnesses can be described in computational terms. Witnesses are given as a black box or an oracle and need not be constructive. When given a witness, you may give it input and will receive appropriate output; you will get no other information.  Multiple instances of a witness can be run, and any instance can be copied.  When requested a witness (as input to a witness), you simply accept input and provide appropriate output; you may give false output, but that may cause the instance of the witness to enter an infinite loop or give arbitrary output. A witness may take arbitrarily long to give an answer and may call appropriate witnesses (for example, a witness to $A\rightarrow B$ may call a witness to $A$) an unlimited number of times.  A witness for a statement containing `$\rightarrow$' can also be viewed as a functional of a higher type.

By running multiple instances of a witness $W$ for $A$, we obtain the response tree for $W$, from which we can obtain $W$ as a real number.  Note that $W$ has to meet the requirements for every path through the tree.  (The response tree specifies the behavior for every possible input.)  For a different notion, see Section 5 Narrow Constructive Truth.

Constructive truth for intuitionistic analysis is defined in Section 6.  Recursive analysis can be interpreted in arithmetic and is not considered separately in this paper.

\section{Basic Properties}
\subsection{Basic Properties}
Some statements, such as $\neg A$, have no existential components (the witnesses have nothing non-trivial to do) and are thus treated classically: $\neg A \leftrightarrow \Box \neg A$ and $\neg A \leftrightarrow \neg \neg \neg A$. $\neg A$ may be viewed as an abbreviation of $A \rightarrow  0=1$. $A\vee B$ is constructively true iff $A$ is constructively true or $B$ is constructively true; analogously with $A\wedge B$ and (for integers) with $\exists x A(x)$. $\forall x A(x)$ is constructively true iff for every $n$, $A(n)$ is constructively true and uniformly so, that is the witnesses can be constructed uniformly in $n$. 

The requirement on witnesses for implication is precisely such that key properties of implication are constructively true. $A\rightarrow B$ is constructively weaker than $\neg A \vee  B$, which allows $A\rightarrow A$ to be constructively valid. If $A \wedge  (A\rightarrow B)$ is constructively true, then so is $B$. $A \leftrightarrow B$ means $(A \rightarrow B) \wedge (B \rightarrow A)$,  and `$\leftrightarrow$' has the same precedence as `$\rightarrow$'.  As required, constructive truth of $A\rightarrow B$ is preserved if the notion of constructiveness is broadened, which is why the witnesses for antecedents cannot be required to be constructive.

Because the whitespace in witnesses can absorb time complexity, every constructively true statement technically has a polynomial time computable witness.  The transformations on witnesses corresponding to the connectives and quantifiers are also recursive (and hence polynomial-time computable).  For example, to convert a witness for $\forall x A(x)$ into a witness for $A(n)$, pick all of the input-output pairs starting with $n$, and remove $n$ from the inputs.

Like in the classical logic, every occurrence of a subformula is either positive or negative. If $A\rightarrow B$ is constructively valid, then any positive occurrence of $A$ may be replaced $B$, and any negative occurrence of $B$ may be replaced by $A$, and the above is constructively true.

Intuitionistic predicate logic, Heyting arithmetic, and Markov's principle $\forall x(A(x)\vee \neg A(x)) \wedge  \neg \forall x\neg A(x) \rightarrow  \exists x A(x)$ are constructively true.  $\forall x(\neg \neg A(x)\vee \neg A(x))$ is constructively true iff $A$ is decidable.

Similarly, $\Box(\forall x(\neg \neg B(x)\vee \neg B(x)) \rightarrow \forall x(\neg \neg A(x)\vee \neg A(x)))$ means $A$ is recursive in $B$.

Every statement $A$ has the associated set of witnesses $S(A)$. $A\rightarrow B$ is constructively true iff there is a recursive mapping that is total on $S(A)$ and maps $S(A)$ into $S(B)$.

Conjunction and disjunction are not strictly necessary since with some modifications to the formula, they may be replaced by the universal and existential quantifier, respectively. 

The property of being a constructive witness is $\Pi^0_2$ above the (constructive) complexity of blocks that start with negation or top-level implication. Thus, constructive truth is $\Sigma^0_3$ above the complexity of such blocks.  Moreover, $\Box \forall x \exists y \phi (x, y) \leftrightarrow \Box \exists n \forall x (\lbrace n \rbrace (x)$ exists $\wedge \: \phi(x, \lbrace n\rbrace(x)))$ is constructively valid, and analogously with more quantifiers.  Here $\lbrace n \rbrace (x)$ is the output of the $n$th Turing machine with input $x$ when that machine halts.

Also, a (constructively) $\Sigma^0_3$ statement is true iff it is constructively true. In fact, $\phi \leftrightarrow \Box \phi$ is constructively valid for (constructively) $\Sigma^0_3$ $\phi$.

Recursive realizability, or just realizability, is like constructive truth except that witnesses must be given as standard codes for partial recursive functions. A statement is realizable iff it has a witness in this sense. The realizability interpretation corresponds to inserting `$\Box$' at every relevant place, that is in front of the statement and antecedents and after `$\neg$'.  For every statement $A$, $A$ or $\neg A$ is realizable (as is $A \vee \neg A$).  The set of realizable sentences is a consistent completion of some constructive formal systems; it has the same Turing degree as the truth predicate for arithmetic.  Intuitionistic predicate logic, Heyting arithmetic, and Markov's principle are realizable.  Some realizable statements are false; an example is ``not every partial recursive function is either defined or undefined at zero''. Extended Church Thesis is the realizable assertion (actually, a schema) that truth equals realizability.

A constructive proof provides a number and demonstrates that the number codes a constructive witness.  (However, even when the constructive proof is short, an actual computation of say which disjunct is true may be unfeasibly long.) A true constructive formal system is a formal system in which only constructively true statements are derivable, and in which derivations can be constructively converted into witnesses for constructive truth. When one is more interested in the witnesses than in the truth, using a constructive formal system may be a good way to manage complexity. For example, constructive formal systems have applications in computer science since algorithms can be ``read off'' constructive proofs. When one is more interested in realizability than in truth, a constructive formal system that proves false statements but does not prove unrealizable statements may be useful. However, the savings have to be balanced against the difficulty of not using some logical truths; and to avoid contradiction, one must always note when ``A'' is used as a shorthand for ``A is realizable''.

\subsection{An Abstract Theory of Constructive Truth}
Constructive truth acts as a necessity operator in intuitionistic modal logic.  If the language is treated constructively, then $A \Rightarrow \Box A$ is admissible as a rule of inference.  We define the abstract theory of `$\Box$' to consist of intuitionistic predicate logic, rule  $A \Rightarrow \Box A$, and axioms $\Box (A\rightarrow B)\rightarrow (\Box A\rightarrow \Box B)$, \enspace $\Box A\rightarrow A$, \enspace $\Box A\rightarrow \Box \Box A$, \enspace $\neg A\rightarrow \Box \neg A$, \enspace $\Box (A\vee B) \rightarrow \Box A\vee \Box B$, and $\Box \exists x A(x) \rightarrow  \exists x \Box A(x)$. Constructive arithmetical truth satisfies the theory.  In the above, we have assumed that every element is constructive (and given as a construction); see Section 6 for a treatment of non-constructive elements.  
The theory proves $\Box (A\wedge B) \leftrightarrow \Box A\wedge \Box B$, \enspace $\Box (A\vee B) \leftrightarrow \Box A\vee \Box B$, \enspace $\Box \exists x A(x) \leftrightarrow \exists x \Box A(x)$, and $\Box \forall x A(x) \rightarrow \forall x \Box A(x)$; but $\neg \neg \Box A \rightarrow A$ need not be constructively valid, that is the universal closure need not be constructively true.

The theory seems to capture the essence of the constructive truth operator relative to the underlying logic.  However, intuitionistic logic has natural extensions.  If the interpretation of `$\forall$' is compatible with the classical one (as is the case here), then $\forall x\neg \neg A(x) \rightarrow \neg \neg \forall xA(x)$ can be added.  If the domain is constructively enumerable and unbounded searches can be performed (as is also the case with constructive arithmetical truth), then $\forall x(A(x) \vee B(x)) \wedge \neg \forall x B(x) \rightarrow \exists x A(x)$ can be added.  However, this principle can fail if the domain consists of say codes for total recursive functions or if we impose time limits on computations.  Other principles can be added, but they should be constructive so that the rule $A \Rightarrow \Box A$ can be retained.

An abstract theory of constructive knowledge consists of the abstract theory of `$\Box$' except for $\neg A\rightarrow \Box \neg A$ but with additional axioms as appropriate.  $\Box A$ would stand for ``$A$ is constructively knowable''.  One notion of constructive knowledge is that a witness for $\Box_P A$ codes a witness for $\Box A$ and a proof (in a sufficiently strong formal system) that the code witnesses $\Box A$.  Another interesting notion of constructive knowledge includes $\neg\neg\Box A \rightarrow A$ (and hence $\Box \forall x(\neg\neg\Box A(x) \rightarrow \Box A(x))$).

\section{Complexity of Constructive Truth}
A $\Pi^1_1$ formula is one that can be presented in the form $\forall X \phi$  where $\phi$  is arithmetical, allowing the use of the set of integers $X$ as a predicate. Thus, a statement is $\Pi^1_1$ if it is equivalent to a statement that a particular recursive (or arithmetical, or polynomial time computable) partial order (or total order) is well-founded.

\begin{description}
\item[Theorem 1:] Constructive arithmetical truth is $\Pi^1_1$ hard.
\item[Proof:]
Let $x = (x_0, x_1, x_2, \ldots)$ be arithmetical and such that for every $n$, $x_n$ (which is a set of natural numbers) is not recursive in $(x_0, x_1, x_2, \ldots)$ with $x_n$ replaced by $0$.  Such $x$ exists by standard results in recursion theory.  For example, every real number Cohen generic over recursive sets can be split into infinitely many mutually generic predicates $x_n$.  

Let $P$ be $\forall i (\neg \neg x_n(i) \vee  \neg x_n(i))$.  We have $\neg \Box \exists n (\forall m\neq n P(m) \rightarrow  P(n))$.

Let $T$ be a recursive tree of integer sequences under extension.  We claim that $T$ is well-founded iff $\Box (\forall n \exists s ($sequence $s$ has length $n \wedge s \in T \wedge P(s)) \rightarrow \exists s \exists t (s$ is incompatible with $t \wedge P(s) \wedge P(t)))$.    If $T$ is well-founded, then we simply wait until incompatible $s$ and $t$ are given, and then copy the witnesses for $P(s)$ and $P(t)$.  If $T$ is not well-founded, then a ``bad'' witness for the antecedent will start to provide an infinite path through $T$ and then once we commit to incompatible $s$ and $t$, provide witnesses for every $P$ except for $P(s)$ or for $P(t)$.
\end{description}

The proof shows that constructive truth for arithmetical statements of the form $A \rightarrow B$ where $A$ and $B$ are without `$\rightarrow$' is $\Pi^1_1$ complete.  Also, by Proposition 3 (below), constructive truth for `$\exists$'-free statements is $\Pi^1_1$ hard.  


Transfinite induction for P on `$\prec$', TI(P, $\prec$) asserts $\forall n(\forall m\prec n P(m) \rightarrow  P(n)) \rightarrow  \forall n P(n)$.

\begin{description}
\item[Proposition 2:]  For every formula P and well-founded relation `$\prec$', $\Box$TI(P, $\prec$).
\item[Proof:]  To get a witness for $\forall n P(n)$, request a witness for $\forall m\prec n P(m) \rightarrow  P(n)$ for every $n$. Then, fulfill all requests for witnesses for $P(m)$ as witnesses become available. By induction on the rank of $n$, for every $n$, a witness for $P(n)$ will be provided.  Witnesses for $m \prec n$ are provided in $\forall m (m\prec n \rightarrow  P(m))$, so `$\prec$' need not be recursive.
\end{description}

Conversely, a form of constructive induction implies well-foundness.

\begin{description}
\item[Proposition 3:]  If `$\prec$' is a recursive strict partial-order given as such, and $R$ is a universal $\Pi^0_1$ predicate, then ``$\prec$ is well-founded''$\leftrightarrow \Box \forall n TI(R(n) \vee \neg R(n), \prec)$.
\item[Proof:]  An input-output pair for $TI(R(n)\vee \neg R(n), \prec)$ can be converted into a finite boolean combination of $\lbrace R(n, 0), R(n, 1), R(n, 2), \ldots \rbrace$ asserting ``input implies output''.  For a witness $S$, every such combination included must be true.  We show that if `$\prec$' is ill-founded, then some combination included is not tautological.  Let $a_0 \succ  a_1 \succ  a_2 \succ \ldots$ be infinite and consider the following witness $T$ to the antecedent:  For every $i$, require the correct truth-value of $R(n, a_{i+1})$ before giving away $R(n, a_i)$; other than that, answer every $R(n, m)$.  Some (segment of $T$, ($a_0$, truth value of $\neg R(n, a_0)$) is included in $S$, but any such combination is non-tautological.  Now, from any witness to $\forall n TI(R(n) \vee \neg R(n), \prec)$, a non-tautological combination of $(R(n, i))_i$ can be extracted for every $n$, so by universality of $R$, the witness is non-recursive, which completes the proof.
\end{description}

In essence, the proof of Proposition 3 uses quantification over all $\Pi^0_1$ predicates to approximate the idealized unknowable predicate $G$:  A constructive witness for a statement using `$G$' must work regardless of which set is $\lbrace n:  G(n) \rbrace$.  The same idea is used to prove the next theorem.

\begin{description}
\item[Theorem 4:] Constructive arithmetical truth is $\Sigma^1_1$ hard.
\item[Proof:] Let $T$ be a recursive tree of integer sequences under extension.  Let $R$ be a universal $\Pi^0_1$ predicate, and $P(a, b)$ be $R(a, b) \vee \neg R(a, b)$.  We prove that $T$ is ill-founded iff:
$\Box \forall n ((\forall l \exists s ($sequence $s$ has length $l \wedge s \in T \wedge P(n, s)) \rightarrow \exists s \exists t (s$ is incompatible with $t \wedge P(n, s) \wedge P(n, t))) \rightarrow S(n))$, where $S(n)$ asserts ``there is a non-tautological finite boolean combination of $\lbrace R(n, 0), R(n, 1), R(n, 2), \ldots \rbrace$ that is not false.''  By universality of $R$, we have $\neg \Box \forall n S(n)$.  If $T$ is well-founded, then the antecedent has a recursive witness, and hence the full statement is not constructively true.  On the other hand, if $T$ is not well-founded, then any witness for the antecedent will have to respond to some branch $s$ of the tree, and hence divulge a statement of the form $\forall t (t$ is a subsequence of $s) (\neg)R(n, t) \rightarrow (\neg)R(n, u)$, where $u$ is not a subsequence of $s$.  Here `$(\neg)$' means that `$\neg$' can independently be either present or absent (but the statement must specify when `$\neg$' is present).  Now, that statement is not a tautology and hence witnesses $S(n)$.
\end{description}

In contrast to theorems 1 and 4, we have the following theorem, which shows that in terms of the number of implications, the complexity of the formulas used in the proofs of the theorems above is optimal.

\begin{description}
\item[Theorem 5:] Let $\phi$ be arithmetical and such that for every $A \rightarrow B$ in $\phi$, $A$ uses neither `$\rightarrow$' nor `$\exists$', but allowing bounded `$\exists$'.  (Negations do not count as implications.)  Then, $\Box \phi$ is constructively equivalent to an arithmetical formula.
\item[Proof Sketch:]  $\Box \phi \leftrightarrow \exists n$ ``$n$ is a witness to $\Box \phi$''.  For statements without `$\rightarrow$', strengthen the notion of witnesses by requiring that answers be given in order and without whitespace.  Any witness can be so strengthened, so constructive truth is unaffected, and we show that under this change, being a witness to $\phi$ is arithmetical.  $W$ witnesses $\phi$ iff ($W$ has the right syntactic form and) for every numerical input into $\phi$ for all witnesses input into $\phi$, appropriate output is given.  $X$ witnesses $A$ (in the strengthened sense) where $A$ is an antecedent inside $\phi$ iff every initial segment of $X$ is correct.  By compactness, $W$ witnesses $\phi$ iff for every numerical input into $\phi$, there is a finite tree of (finite) witness inputs into $\phi$ such that $W$ gives appropriate and complete numerical output for every witness input outside the finite tree (allowing false or incomplete output for false input).  The above formula is arithmetical, which completes the proof.
\item[Corollary to the Proof:] Constructive truth for arithmetical statements for which every antecedent $\phi$ is as $\phi$ is in the theorem is $\Pi^1_1$.
\end{description}

We conjecture constructive arithmetical truth to be complete for second order arithmetic and have the following partial results:

\begin{description}
\item[Proposition 6:]  The relation of being a witness is complete for second order arithmetic.
\item[Proof:]   By inspection, for every $\phi$, being a witness to $\phi$ is definable in second order arithmetic.  Let $A_0$ be $\forall x \exists y \: 0=0$, and $A_{n+1}$ be $A_n \rightarrow (0=0 \vee 0=0)$.  (By coding in whitespace, we could even have started with $A_{-1}$ as $0=0 \vee 0=0$, or without using whitespace, $A_{-1}$ as $\forall x (0=0 \vee 0=0)$.)  We define a recursive function $f$ such that for every infinite binary sequence $X$ and every $\phi$ given as $\Pi^1_n$ ($n>0$) but not as $\Pi^1_{n-1}$, $\phi(X)$ holds iff $f$(``$\phi$'', X) is a witness to $A_n$.  If $\phi$ is $\Pi^1_1$, convert it to an assertion that a particular recursive-in-$X$ relation is well-founded.  Use the witness for the antecedent to guess an infinite descending path, and when that fails, give a witness to the consequent.  Also, use additional delays to code $X$ into $f$(``$\phi$'', X).  If $\phi$  is $\Pi^1_{n+1}$ (but not $\Pi^1_n$ or $\Pi^1_1$), then convert it to $\forall Y \psi(X$ join $Y)$ where $\psi$ is $\Sigma ^1_n$, and delay giving a witness to the consequent for as long as the candidate witness for the antecedent is $f$(``$\neg \psi$'', X join Y).
\end{description}

\begin{description}
\item[Theorem 7:] Constructive arithmetical truth is is $\Pi^1_2$ and $\Sigma^1_2$ hard.
\item[Proof:] We build on the proofs of Theorem 4 and Proposition 6.  As in Theorem 4, we quantify over all $\Pi^0_1$ predicates.  Let `$R$' be a symbol for a $\Pi^0_1$ predicate.  For simplicity, we present our statements using $R$, and interpret $\Box D$ as $\Box \forall n D(R'_n)$ where $R'$ is a universal $\Pi^0_1$ predicate.  Let $P(n)$ be $R(n) \vee   \neg R(n)$.  The $\Pi^1_2$ hard formula $\Box D$ will have $D = A\wedge B\rightarrow C$. $A$ will assert $\forall n\exists s ($sequence $s$ has length $n \wedge s\in T \wedge P(s)) \rightarrow  \exists s\exists t (s$ is incompatible with $t \wedge P(s) \wedge P(t))$, where $T$ is the tree of integer sequences that is effectively coded (say, as a primitive recursive predicate) by the free variable of $A$, and with sequences coded by odd natural numbers.  $B$ asserts $\forall n (P(4n) \vee P(4n+2))$.  $B$ is used as a second order universal quantifier; witnesses for $B$ code real numbers through choices of $4n$ or $4n+2$.  $C$ states $\exists s \exists u (\forall n (n$ in sequence $s) P(2n) \: \wedge $ sequence $u$ is not a subsequence of $s \: \wedge (\forall t (t$ is a subsequence of $s) (\neg )R(t) \rightarrow  (\neg )R(u)))$.  

We claim that $\Box D$ holds iff for every witness to $B$, there is an infinite path through $T$ such that for every $n$ in the path, $P(2n)$ is given by the witness.  To see that the claim implies the theorem, note that the statement is $\Pi^1_2$ complete (as a function of the parameter of $A$).  For $\Sigma^1_2$ hardness, let $E$ be $D \rightarrow$ ``there is a true non-tautological finite boolean combination of $\lbrace R(0), R(1), R(2), \ldots \rbrace$''.  Clearly, by universality of $R$, $\Box E\rightarrow \neg \Box D$, and in the proof of the claim, we will show $\neg \Box D\rightarrow \Box E$.

Now, consider a witness $X$ to $B$. Every witness to $A$ will have to respond to every infinite path through $T$ corresponding to $X$.   If the $\Pi^1_2$ statement is true, then such a path exists and hence (like in the proof of Theorem 4), by reading the witness for $A$, one can extract the information required for $C$.  If the $\Pi^1_2$ statement is false, let $X$ be a counter-example witness for $B$. A ``bad'' witness $Y$ for $A$ will require correct answers to $P$ and will delay giving away incompatible elements until the path used as input into $A$ is past the point where some $P(2n)$ is not given by $X$, and then give out $P(s)$ and $P(t)$ for incompatible $s$ and $t$ that are independent of the path.  For every witness $Z$ to $D$, the input-output pair corresponding to $X$ and $Y$ is not tautological. Thus, a non-tautological input-output pair can be extracted from $Z$ through complete search, which completes the proof.

\item[Question:]  What is the complexity of constructive arithmetical truth?
\end{description}

We expect that the proof of Theorem 7 can be generalized to prove that constructive arithmetical truth is complete for second order arithmetic.  For $\Pi^1_{n+1}$ hardness, let $A_{n+1}=(A_n \wedge B_n \rightarrow C_n)$ with $A_1=A$, $B_n$ analogous to $B$ (but using a different domain for each $n$), and $C_n$ asserting existence of a true non-tautological combination of a certain type. However, there are combinatorial difficulties in actually stating a correct $C_n$, so we do not know how to prove the conjecture.

We note that $A_n$ use the optimal number of implications since constructive truth for statements with implication nesting depth $n>0$ (or depth $n+1$ for statements corresponding to Theorem 5) is $\Pi^1_n$.

\begin{description}
\item[Question:] Are there arithmetical statements not constructively equivalent to `$\exists$'-free ones?
\end{description}

Probably, the answer is yes and $\exists n P(n)$ is such a statement where $P$ is as in Theorem 1. On the other hand, we expect that every $\Box \phi$ is constructively equivalent to $\Box \psi$ where $\psi$ uses neither `$\exists$' nor `$\Box$'.  If $\Box A(n)$ is equivalent to ``$n$ is a witness to $\phi$'', then $\Box \phi$ is equivalent to $\Box \exists n A(n)$ and may be equivalent to $\Box (\forall n(A(n) \rightarrow B) \rightarrow B)$ for appropriate $B$. $A$ (and $B$) can likely be chosen without `$\Box$' and `$\exists$'.

\section{Narrow Constructive Truth}
Narrow constructive truth is a particular modification of constructive truth in which witnesses are not completely given.  A witness for narrow constructive truth provides a code for a recursive interactive function such that for every possible interaction, the witness requirements are met.  The meeting of requirements is defined inductively on the complexity of the formula.  If $A$ is false, then the requirements for $A$ are immediately violated.  For $A\rightarrow B$, multiple numbered instances can be launched, and different choices for the candidate witness for $A$ require different instances.  A new instance may copy the input and output (but not the internal state) from any point of an existing instance (it is unclear if omitting this provision would change which statements are narrow constructively true).  Every instance for $A\rightarrow B$ must meet the requirements for $B$ for as long as the requirements for $A$ are met, including in the limit (that is at time infinity).  If every instance meets the requirements, then the requirements for $A\rightarrow B$ are met even if answers are omitted along some path where the copying was done infinitely many times.  For other connectives, the meeting of requirements is as usual.  For example, in $A \vee B$, one must choose $A$ or $B$, and then meet the requirements for the choice.  Failure to choose violates the requirements at infinity but not at any finite time.

Narrow constructive truth is similar to constructive truth, with the difference limited to statements with implications inside antecedents.  It satisfies all of the basic properties listed (except for the description of $A \rightarrow B$ in terms of mappings), and theorems 1, 2, 3, and 5.  Narrow constructive truth implies constructive truth for statements without antecedents inside antecedents inside antecedents, and perhaps for all arithmetical statements.  \textit{For this section only}, `$\Box$' represents narrow constructive truth.

A key distinguishing feature of narrow constructive truth is that for $A$ and $C$ without `$\rightarrow$', $\Box((A\rightarrow B)\rightarrow C) \wedge\neg\Box A\wedge\neg\neg B \rightarrow \Box C$ is narrow constructively valid (equivalently, constructively valid).

\begin{description}
\item[Theorem 8:] Narrow constructive arithmetical truth is $\Pi^1_1$ complete.
\item[Proof:]
By the proof of Theorem 1, narrow constructive truth is $\Pi^1_1$ hard.  To see that it (and the corresponding notion of witnesses) is $\Pi^1_1$, a statement is narrow constructively true iff there is an integer (coding the narrow constructive witness) such that for every real number interacting with and trying to refute the alleged constructive witness, the witness requirements are met.   The meeting of requirements is hyperarithmetically definable.
\end{description}

Thus, by using a universal $\Pi^1_1$ formula, we can find a ``universal'' arithmetical formula $\phi$  such that for every arithmetical sentence $\psi$, $\Box \phi(`\psi$'$)$ iff $\Box \psi$. This does not contradict undefinability of truth since for some `$\psi$', $\phi(`\psi$'$)\rightarrow \psi$  is false.

The reasoning can be used to compute the expressive power of the full language, which allows nesting of the narrow constructive truth operator. A hyperjump of $X$ is essentially the set of (standard codes for) well-founded relations that are recursive in $X$. For every natural number $n$, there is a formula $\phi$  that only uses `$\Box$' $n$ times, such that $\Box \forall x(\Box \phi (`\psi (x)$'$)\leftrightarrow \Box \psi (x))$ (and hence $\Box \phi (`\psi$'$)\leftrightarrow \Box \psi$  if `$\psi$' is a statement) whenever `$\psi$' has at most $n$ nested `$\Box$'; and $\{x:\Box \phi (x)\}$ has the same complexity (even under one-to-one polynomial time reducibility) as the (n+1)th hyperjump of 0. The proof is by induction on n, and is analogous to the proof of Theorem 8. The proof also shows that for expressible $P$ and $\phi$, ``$P$ is a witness to $\phi$'' is also expressible.

\section{Constructive Truth for Analysis}
\subsection{Prelude and Definition}
Before defining constructive analysis, we note that analysis can be interpreted in second order arithmetic.  Analysis deals not with arbitrary functions but continuous and other reasonable functions, and these can be represented by real numbers.  Under an appropriate representation, the theorems of analysis are almost all $\Pi^1_3$, and are usually $\Pi^1_2$.  Moreover, most of the $\Pi^1_2$ theorems have a natural $\Pi^1_1$ strengthening by requiring that in $\forall X\exists Y \phi(X, Y)$ $Y$ is arithmetical in $X$.  Similarly, most $\Pi^1_3$ theorems $\forall X\exists Y \forall Z \phi(X, Y, Z)$ in analysis have a $\Pi^1_2$ strengthening by requiring $Y$ to be recursive in a finite hyperjump of $X$.  Thus, analysis can essentially be interpreted through constructive truth in arithmetic.  Philosophically, one can argue that natural numbers and constructive arithmetical truth are basic, so the high complexity of constructive arithmetical truth provides an (additional) ontological basis for analysis.

We use the language of second order arithmetic.  Capital letters represent sets of natural numbers.   Free variables need not be shown in schemas.

While there are different notions of constructive analysis, prominent among them is Brouwer's intuitionistic analysis.  Under it, infinite binary sequences are given not as completed objects, but as a process that gives out zeroes and ones.  Under the concept of choice sequences, we can make completely arbitrary choices as which natural numbers belong to a choice sequence $X$.  We construe this to mean that any real number can be a choice sequence.  Thus, the domain will include all sets of natural numbers (see Section 7 for discussion).  Other approaches treat choice sequences either as a figure of speech or as a fundamentally non-classical construct.  Also since the sets represent real numbers rather than being a shorthand for arbitrary properties, we will have $\forall X \forall n (X(n) \vee \neg X(n))$.

Answers must be given in finite time rather than after the process is completed.  Therefore, we reinterpret $\forall X$ to mean for all $X$ and with answers given continuously in $X$.  Under this view, every function is continuous on its domain, and $\forall X(X=0 \vee  \neg X=0)$ is false since one cannot always decide whether $X=0$ by looking at only a finite segment of $X$.  Since the statements are not treated classically, this is not a contradiction.  We use $W(\phi)$, $\phi$ is witnessable, for this interpretation of $\phi$.  We define ``$\phi$ is continuous function realizable'' to mean $W(\phi)$.

We formally define witnesses like in Section 2, except that being witnessable need not match classical truth.  However, we do have $W(\neg \phi)\leftrightarrow \neg W(\phi)$.  Input in input-output pairs includes initial segments of the sequences for the second order universal quantifiers, and output includes initial segments of the sequences for the existential quantifiers.  $W$ witnesses $\forall X \phi (X)$ iff for all $X$, $W_X$ witnesses $\phi(X)$.  As in Section 2, $W_X$ is obtained from $W$ be selecting input-output pairs corresponding to $X$ and and removing $X$ from input.  $W$ witnesses $\exists X \phi(X)$ iff $W$ provides arbitrarily long initial segments for a single $X$, and a witness for $\phi(X)$, but not initial segments incompatible with $X$.  Treatment of other connectives (except `$\Box$') is unchanged.

Existence of non-constructive reals requires a change in the treatment of `$\Box$'. A witness for $\Box_S \phi$ provides a standard code for a recursive in $S$ witness to CL$_S (\phi)$ where CL$_S$ is the universal closure except for first order variables and for variables in $S$.  $S$ is written as a list of second order variable symbols and is omitted when empty.  $S$ is treated as a set; for example, $\Box_{X,Y} \phi \leftrightarrow \Box_{Y,Y,X} \phi$.  Also, the semantics of `$\Box$' is definable in second order arithmetic, so we can treat `$\Box$' not as a primitive symbol but as a particular transformation of formulas.

\subsection{Basic Properties} 
To define an abstract theory of $\Box_S \phi$, we can add `$\leq_c$' with $S\leq_c T$ meaning $S$ can be constructed from $T$.  We allow $\leq_c$ to act on finite sets of elements, but with each element listed. We can axiomatize the theory by using the abstract theory of `$\Box$' for every $\Box_S$ (except that in $\neg A \rightarrow \Box_S \neg A$, we require $S$ to include all free second order variables of $A$), properties of $\leq_c$ (reflexivity, transitivity, $S,T$ being the $\leq_c$-least upper bound of $S$ and $T$, and the empty set being the $\leq_c$-least element (up to $\leq_c$ equivalence)), $\Box_S\phi \rightarrow \Box_S \forall X\phi$  ($X$ not in $S$), $S\leq_c T \rightarrow (\Box_S\phi \rightarrow \Box_T\phi)$ (no element of $S\backslash T$ is free in $\phi$), and $\Box_S \exists X \phi \rightarrow \exists X \leq_c S$ $\Box_{S,X} \phi$  ($X$ not in $S$).  Our constructive analysis satisfies additional properties including pairing and $\phi \rightarrow \exists X\Box_{S,X}\phi$ (completeness for witnessability; here $S$ includes all free second order variables of $\phi$), and in it, $S\leq_c T$ means $S$ is recursive in $T$.  

Because of continuity (which is classically false) and compactness of $2^N$, we even have $\Box(\Box_S \forall X \phi \leftrightarrow \forall X \Box_{S,X} \phi)$ ($X$ not in $S$).

The entire discussion of Sections 3-5 generalizes to formulas with an arbitrary second order parameter $X$ by replacing recursive with recursive in $X$ and analogously replacing $\Pi^0_1$ and other notions, and using $\Box_X$ in place of $\Box$.

We note that a witness for $\forall X (\phi (X)\rightarrow \psi (X))$ is simply a code for continuous function $f$ such that if $Y$ is a witness to $\phi(X)$, then $f(X, Y)$ is a witness to $\psi(X)$.  (Because we allow whitespace in witnesses, we can require $f$ to be total.)

Formulas without `$\exists$' and `$\vee$' (and, but this is unnecessary here, with all atomic formulas negated) are called negative formulas and are treated classically.  Every formula is classically equivalent to a negative one, so classical analysis can be interpreted in the constructive one.  Conversely, the relation of being a witness is definable in analysis, so constructive analysis can be interpreted in the classical one.  

Constructively, we define $W(\phi)$ as $\exists X$ ``$X$ is a witness to $\phi$'', where the quoted statement is negative.  We have $\Box (\phi\leftrightarrow W(\phi))$, so constructively, every formula is equivalent $\exists X$ ``a negative formula'', and every negation equivalent to a negative formula.  We expect that every $\Box_S \phi$ is constructively equivalent to $\Box_S \psi$ for an arithmetical $\psi$ without `$\exists$'.

Because the relation of being a witness is definable, we can define additional connectives by specifying witness requirements in terms of witnesses for the constituent formulas.  We define a weak or non-constructive `or' by:  A witness for $A\vee_w B$ must provide $X$ and $Y$ such that $X$ witnesses $A$ or $Y$ witnesses $B$.  The following are constructively valid:  $A\vee_w B \leftrightarrow B\vee_w A$, \quad $A\vee_w (B\vee_w C) \leftrightarrow (A\vee_w B)\vee_w C$, \quad $\neg A\vee_w \neg\neg A$, \quad $A\vee B \rightarrow A\vee_w B$, \quad $(A\vee_w B) \wedge \neg B \rightarrow A$, \quad $(B\rightarrow C)\rightarrow (A\vee_w B\rightarrow A\vee_w C)$, \quad $\forall X(\neg A(X)\rightarrow B(X))\vee_w C \leftrightarrow \forall X(\neg A(X)\rightarrow B(X)\vee_w C)$ (and hence $A \wedge B \vee_w C \leftrightarrow (A\vee_w C) \wedge (B\vee_w C)$), \quad $\Box_S(A\vee_w B) \rightarrow \Box_S A \vee_w \Box_S B$, and (because of our particular treatment of witnesses) $(A\vee_w B \rightarrow \Box_S C) \wedge \neg \neg (A \wedge B) \rightarrow \Box_S C$, but not $\Box(A\vee_w A) \rightarrow A$ nor $\Box A \vee_w \Box A \rightarrow \Box(A\vee_w A)$.  `$\vee_w$' is the strongest (definable) connective satisfying $\Box ((\neg A \rightarrow B) \wedge (\neg\neg A \rightarrow C) \rightarrow B\vee_w C)$.  

A generalization of `$\vee_w$' is given by:  A witness for $\exists_w n A(n)$ must provide $X_0, X_1, X_2, \ldots$ such that for some $n$, $X_n$ witnesses $A(n)$.  The quantifier $\exists_w$ satisfies the analogues of all of the above properties of `$\vee_w$', and also $\Box_S \exists_w n A \leftrightarrow \neg\neg\Box_S A$ ($n$ not free in $A$).  For example, $\forall n (A\rightarrow B)\rightarrow (\exists_w n A \rightarrow \exists_w n B)$ is constructively valid, which means that $\exists_w$ is non-modal and monotonic.  A formula $\phi(..)$ is non-modal iff for all $A$ and $B$, $Cl_\phi (A\leftrightarrow B)\rightarrow (\phi(A) \leftrightarrow \phi(B))$ is constructively valid, where $Cl_\phi (C)$ is the universal closure for variables free in $C$ that are not free in $\phi(C)$.

We can also formalize the unknowable predicate $G$ by treating $\forall G \phi(G)$ as $\forall X \phi(G)$ where $G(m)$ is ``$\lbrace n: (m, n) \in X \rbrace$ is not finite'' and X is not otherwise used in $\phi$.

A constructive system for analysis can be defined as one in which proofs of $\phi$ can be constructively (and uniformly in $\phi$) converted into (constructive) witnesses for the universal closure of $\phi$.  Continuous function realizability can also be treated classically, but then one would not have $W(\phi)$ to $W(\forall X\phi)$ as a rule of inference.

If $\phi$ is false, then $W(\phi)$ can only hold if $\phi$ has a negative occurrence of $\forall X A(X)$ with $A$ containing `$\exists$' or `$\vee$' not in the scope of an antecedent (or negation) inside $A$ ($X$ can be any second order variable).  Thus, we can form a classically correct subsystem by including all derivable (or witnessable) formulas not of this type, and using classically correct principles for other formulas.  In constructivizing theorems of analysis, there is usually no problem in avoiding that type of formulas.  Principles that are both classically and constructively correct include induction, transfinite induction, countable choice, and dependent choice.  Apparently, even adding $\neg\neg\forall X(A(X)\vee\neg A(X))$ will result (even in the presence of the above principles) in a system that is constructive with respect to formulas $\phi$ in which every positive (in $\phi$) occurrence of $\forall X$ that is inside an antecedent or negation is part of `$\forall X \neg$'.

\subsection{Development of Constructive Analysis}

A key constructively valid principle is countable choice $\forall n\exists X \phi(X, n)\rightarrow \exists X\forall n\phi(X_n, n)$.   Countable choice is a rather strong principle and implies constructivized comprehension: $\forall n(\phi(n)\vee \neg \phi(n)) \leftrightarrow \exists X \forall n (n\in X \leftrightarrow  \phi(n))$  ($X$ not free in $\phi$).  In combination with $\neg \neg \forall n(\phi(n)\vee \neg \phi(n))$ (which is also constructively valid), constructivized comprehension is quite powerful and implies the negative form of full comprehension.  Constructively, we also have transfinite induction $\neg \exists$ (an infinite $\prec$-descending path) $\rightarrow $TI$(\phi,\prec)$ for all formulas $\prec$ and $\phi$.  We also have Markov's principle, which simplifies to $\forall X(\neg \forall n\neg X(n) \rightarrow  \exists n X(n))$.

We quantify over real numbers, continuous functions, and infinite sequences of integers (points in the Baire space) by using their effective binary codes and writing for example $\forall X$(``$X$ codes an infinite sequence of integers'' $\rightarrow \phi(X)$). A real number is coded by a nested sequence of intervals with rational endpoints.  By using a conservative and essentially definitional extension, one can also quantify over real numbers (and open sets and continuous functions) directly rather than through coding.  

\begin{sloppypar}Constructively, we have the continuity principle $\forall X\exists Y \phi(X, Y) \leftrightarrow \exists ($continuous $f) \forall X \phi(X, f(X))$, and analogously with more quantifiers.  This also holds for quantification over the Baire space $N^N$.  However, as stated the principle would fail for real numbers and should be modified to use $\exists !Y$ (here, the uniqueness is for the real number coded by $Y$ rather than the code $Y$). This is needed because different sequences can denote the same real number.\end{sloppypar}

In contrast to arithmetic, many basic theorems of analysis are not constructively true (except in the negative form), and most of those that are true have to be stated carefully.  Clearly, constructive analysis should be meant to complement rather than replace classical analysis.  On the other hand, many non-constructive theorems have constructive versions.  For example, the following is a constructive version of the intermediate value theorem: If $f$ is a continuous function (on real numbers) that has both a positive and a negative value, and if, for no interval of positive length, $f$ vanishes on every point of the interval, then $f$ has a zero.

It is crucial that this theorem and other theorems in analysis are stated as implications rather than disjunctions.  One generally cannot decide from the code whether the conditions of the theorem apply.  However, given the conditions, the construction must succeed.  A major reason to study constructive analysis is that it gives effective versions of theorems in analysis.  The added condition for the constructive intermediate value theorem may sound strange, but the point is that under these conditions we can actually compute a zero of $f$ and to do so uniformly in $f$.

Constructively, it is not the case the every bounded monotonic sequence of real numbers has a limit.  To make the theorem constructively true we could include existence of the Turing jump of the sequence (or a certain much weaker property) as a condition.  A crucial substitute that we do have is compactness of the unit interval: Every (countable) open cover has a finite subcover.  Note that the compactness fails in recursive analysis, and the difference is that we require our continuous functions to be continuous even at non-recursive real numbers.  The compactness is a consequence of the Fan Theorem, also called constructive Weak K\"{o}nig's Lemma.  An equivalent of the constructive Weak K\"{o}nig's Lemma is:  Every countable theory that does not have a model (with a complete diagram) has an inconsistency.

Most consequences of the Weak K\"{o}nig's Lemma in analysis rely only on the constructively true form, so in terms of reverse mathematics, constructive analysis corresponds more closely to WKL than to RCA.  For example, even constructively, continuity implies uniform continuity for functions from a complete totally bounded metric space.

\subsection{Constructive Reverse Mathematics}
Constructive reverse mathematics analyzes formal relationships between constructively true statements.  Here, the emphasis is on deriving strong principles in weak theories.

Every constructive theory $S$ has its classical counterpart $T$ consisting of first order logic and negative theorems of $S$.  Conversely, every theory $S$ satisfying ED and containing intuitionistic RCA$_0$ and $\phi \leftrightarrow W(\phi)$ is determined by its classical counterpart $T$.  RCA$_0$ consists of basic arithmetic, $\Sigma^0_1$ induction and recursive comprehension.  ED stands for explicit definability, and this form of ED is sufficient:  If $S$ proves a sentence $\exists X \phi(X)$, then for some $n$, $S$ proves $\exists X$``the $n$th Turing machine outputs $X$ and $\phi(X)$''.  Because any constructive proof can be recursively converted into such a recursive code, any reasonably closed constructive system satisfies ED.  For any classical $T$ containing RCA$_0$, the (for sound $T$, constructive) counterpart $S$ proves a sentence $\phi$ iff for some $n$, $T$ proves ``$n$ witnesses $\Box \phi$''.  S can be axiomatized by intuitionistic RCA$_0$ + $\phi \leftrightarrow W(\phi)$ + negative translation of $T$. $S$ satisfies (the above form of) ED and is conservative over $T$ for negative formulas.  (ED is most often stated to imply $\Sigma^0_1$ soundness for consistent $S$, but the above form of ED allows us to state the general correspondence between $S$ and $T$.) 

Constructive analysis can be done in a weak formal system.  Classical mathematics in weak systems is developed in [Simpson 1999], and doing constructive mathematics is mostly analogous.  Many results can be proved even in intuitionistic RCA$_0$ + $\phi \leftrightarrow W(\phi)$.  The system is strong enough to prove the full continuity principle for Baire space, which implies countable choice and hence constructivized comprehension.  On the other hand, the system is weak enough to be consistent with ``every real number is not not recursive'', an assertion that would turn our analysis into an analogue of constructive recursive analysis except that recursive codes cannot be extracted from real numbers.

\begin{sloppypar}Because of $\Sigma^0_1$ induction, the continuity principle (for Baire space) also implies a form of dependent choice:  $\forall X\exists Y\phi(X,Y) \rightarrow \exists X (X_0=0 \wedge \forall n \phi(X_n, X_{n+1}))$ where $X$ and $Y$ range of the Baire space (this implies countable choice and is classically correct).  Over intuitionistic RCA$_0$ + $\phi \leftrightarrow W(\phi)$, the (constructively) stronger $\exists X\phi(0,X)\wedge\forall X(\exists Z\phi(Z,X) \rightarrow \exists Y\phi(X,Y)) \rightarrow \exists X(X_0=0 \wedge \forall n \phi(X_n, X_{n+1}))$ is equivalent to full induction (equivalently, induction for negative formulas).\end{sloppypar}

Analysis can be developed in intuitionistic RCA$_0$ by using added assumptions about functions such as having a modulus of uniform continuity. Bishop's analysis can essentially be carried out in intuitionistic RCA$_0$ + Countable Choice; and the use of countable choice can likely be eliminated by replacing it with recursive comprehension, and by using $\exists X\forall n\phi(X_n, n)$ rather than $\forall n\exists X \phi(X, n)$ in theorems and definitions (and when appropriate, analogously replacing $\forall m \exists n \phi(m, n)$). 
 
\begin{sloppypar}The schema $\phi \leftrightarrow W(\phi)$ is actually equivalent over intuitionistic RCA$_0$ to Markov's Principle $\forall X(\neg \forall n\neg X(n) \rightarrow  \exists n X(n))$ plus generalized continuity principle: $\forall X(\neg\phi(X) \rightarrow \exists Y\psi(X,Y)) \rightarrow \exists($continuous partial $f) \forall X(\neg\phi(X) \rightarrow $``$f(X)$ is defined'' $\wedge \psi(X,f(X)))$.   The proof of $\phi \leftrightarrow W(\phi)$ is by induction on the complexity of $\phi$.  Markov's principle follows from $\phi \leftrightarrow W(\phi)$ because of the way $W(\phi)$ is formulated. \end{sloppypar}

A classically correct subsystem of intuitionistic RCA$_0$ + $\phi \leftrightarrow W(\phi)$ is intuitionistic RCA$_0$ + Markov's Principle + Countable Choice.  Intuitionistic RCA$_0$ + $\phi \leftrightarrow W(\phi)$ is conservative over the subsystem for formulas in which every negative occurrence of `$\forall X$' classical (that is (provably) replaceable by `$\forall X \neg \neg$').  (The proof is be showing $W(\phi)\rightarrow\phi$ for these formulas.)  Moreover, countable choice is unneeded for formulas with every negative occurrence of `$\forall$' classical. By extension, the conservation result also holds in the presence of additional negative axioms.  This way constructive analysis can be effectively done even in a classically correct subsystem.  Note that while countable choice is classically a rather strong principle, constructively it can be part of a weak formal system.

Most constructive results can actually be proved in just the constructive counterpart of WKL$_0$, that is in intuitionistic RCA$_0$ + $\phi \leftrightarrow W(\phi)$ + constructive Weak K\"{o}nig's Lemma.  The system proves compactness of the unit interval and results that follow from the compactness, including for example, the constructive intermediate value theorem. (That theorem is equivalent over intuitionistic RCA$_0$ to Markov's principle plus constructive Weak K\"{o}nig's Lemma.  In fact, most constructive theorems relying on compactness are equivalent over RCA$_0$ + $\phi \leftrightarrow W(\phi)$ to constructive Weak K\"{o}nig's Lemma.)  Still, the system is weak enough to be $\Pi^0_2$ conservative over primitive recursive arithmetic and to be consistent with ``not every Turing machine either halts or does does not halt on the empty tape''.  

The intuitionistic counterpart of ACA$_0$ is (equivalent to) intuitionistic RCA$_0$ + $\phi \leftrightarrow W(\phi)$ + for arithmetical formulas $\phi$, $\neg \neg \forall n(\phi(n)\vee \neg \phi(n))$.  Finally, over intuitionistic RCA$_0$ + $\phi \leftrightarrow W(\phi)$, the schema $\neg \neg \forall n(\phi(n)\vee \neg \phi(n))$ is equivalent to the negative translation of the full second order arithmetic with countable choice.  For an overview of the reverse mathematics of the law of the excluded middle, see for example [Ishihara 2004]; the base system there can likely be taken as intuitionistic RCA$_0$ + Countable Choice.

We conclude with a principle of constructive reasoning whose (constructive) correctness is not provable in ZFC: $\neg \neg (\forall X\exists Y (\phi(X)\leftrightarrow \psi(Y)) \vee \forall Y\exists X (\phi(X)\leftrightarrow \neg \psi(Y)))$ where $\phi$ and $\psi$ are negative (alternatively, negated) ($\phi$ and $\psi$ may have additional parameters, but all occurrences of $X$ and $Y$ are shown).  Because of the continuity principle, constructive correctness of the schema is equivalent to classical truth of projective Wadge determinacy, which is not provable in ZFC.  Projective Wadge determinacy is implied by and probably equivalent to projective determinacy.

\section{Discussion}

\subsection{General Constructive Mathematics}
Since the law of excluded middle $A \vee \neg A$ is a logical truth, one may wonder in what general sense can it fail.  The answer is that constructively, `$A$' is treated not just as a proposition to be true or false, but also as information or knowledge.  Classically, for every formula and tuple of parameters, the formula either holds for the parameters or not, and no further structure is considered.  Constructively, however, true statements are distinguished according to the information required of them.  Logical connectives, beyond their effect on truth-values, also act on and help specify the information required.  For example, in $\forall n (A(n) \vee \neg A(n))$, the information includes for every $n$ whether $A$ holds.  Constructively, one uses $\neg\neg A$ to assert the truth of $A$ without giving the information corresponding to $A$.

A witness for $A$ consists of all information necessary for $A$; $\neg A$ holds iff $A$ has no witnesses.  Note that the information contained in the witness does not include the reasons that the information is correct; for example, a witness for $\neg A$ just asserts $\neg A$ and need not explain why $A$ fails.  Constructively, a formula with an assignment of free variables is entirely determined by the requirements on witnesses.  The general meaning of logical connectives is as listed in Section 2, and in a way consistent with intuitionistic logic, each framework specifies the details, including what it means to give a witness for the antecedent.  $A$ is constructively true iff a finite amount of knowledge can be transformed (in the required sense) into information sufficient for $A$; and a witness for $\Box A$ supplies that finite information (which can be non-unique).  A constructive justification of $A$ is a classical justification that the integer specified witnesses $\Box A$.

Some constructive frameworks have $\neg \forall n (A(n) \vee \neg A(n))$ for some $A$.  The reason is that the universal quantifier acts not just on the truthfulness of the assertions, but also on the information required.  In these frameworks, the information required for $\forall n (A(n) \vee \neg A(n))$ may (in a specific sense) be impossible, which falsifies $\forall n (A(n) \vee \neg A(n))$.  The classical universal quantifier is `$\forall x\neg\neg$', which ignores which information is required, and `$\forall x$' is a strengthening of that quantifier.   Intuitionistic logic proves $\forall n \neg\neg (A(n) \vee \neg A(n))$, so for formulas with the intensional aspects suppressed (including all negative formulas), intuitionistic logic is simply the full classical logic.

\subsection{Constructive Arithmetic}
The high complexity of constructive truth reminds us of a key philosophical question: Are there other acceptable notions of constructive arithmetical truth? Our answer is no. Other notions such as recursive realizability and narrow constructive truth may be useful, but they do not meet all of the requirements. All of the conditions on the witnesses except the implication one are straightforward. The definition of constructive implication is fixed through a combination of several requirements. Key properties of implication must be true constructively, and constructively true statements must be true. These requirements force us to use non-recursive witnesses in defining constructive truth. This compels us to provide witnesses like real numbers, which essentially fixes how the witnesses are given.

Every true arithmetical statement has an arithmetical witness, but the property of being a witness is not arithmetical. Had we required witnesses to be arithmetical, constructive truth for arithmetical statements would be (mostly) weakened yet remain much stronger than simply truth. However, we cannot do that since we have extended the language of arithmetic,  which requires non-arithmetical witnesses.  We cannot place different complexity requirements (in a truth-altering way) on witnesses for $P$ than for $Q$ since we may discover $P\leftrightarrow Q$ constructively, which can cause inconsistency.

Intuitionistic logic is suitable for various things, and there are different notions of constructive, of which constructive truth is just the default notion.  Another interesting approach to constructivity is to use the full theory of finite types, and treat witnesses for statements with implications as functionals of higher types.  A statement would hold constructively iff it is witnessed by a constructive functional, and constructive functionals would be the recursive ones. By using $\forall x (A(x)\rightarrow \exists y B(x,y)) \leftrightarrow \exists F\forall x (A(x)\rightarrow B(x,F(x)))$ ($A$ and $B$ are negative), every formula would be equivalent to a list of existential quantifiers followed by a negative formula.  This approach would allow making constructive sense of, for example, what is implied by $\forall X(\exists n X(n) \vee \forall n \neg X(n))$ (the ability to test whether any given real number is $0$).  However, for constructivity, we would have to allow functionals to be partial (in which case $B(F(x))$ implies $F(x)$ exists) and either non-extensional or multi-valued, and the semantics is unclear.  Also, in constructive arithmetic, computations should be based on numbers rather than on uncountable objects; so for constructive truth, we have to require that witnesses act based on finite segments rather than on the full witnesses, which leads to a kind of overspill exploited in Theorems 4 and 7.

\subsection{Intuitionistic Analysis}
We believe that the uniqueness extends to analysis.  Our constructive analysis is more usable than any other previously proposed.  (Compare it for example with Troelstra's syntactic treatment of choice sequences.)  Our constructive treatment of implication allows all real numbers to be present in a negative form, yet maintain constructivity.   The naturalness with which our framework extends to analysis confirms its correctness for arithmetic.  

Some critics may question whether our analysis is really constructive.  However, whether or not it represents an intuitionistic foundational stance, our analysis is constructive in that constructive truth of statements amounts to truth of the effective versions of those statements, and in that constructive proof of a statement specifies through a recursive procedure a way that the existential aspects of the statement are satisfied.  For example, existence of a non-recursive real number is not claimed to be constructively true.

Our constructive analysis is derived from two choices, and a reasonable formalizations of the choices.  First, we allow any real number in the domain.  We have to allow every real number that might be generated by a physical or a mental process.  However, our analysis is objective and mathematical (as opposed to being a physical theory), so it has to be independent of physical laws and of human mental capabilities.  Moreover, even under physical laws as currently understood, many physical processes are random, so not every physical process is recursive.  This compels us to include every real number.

Second, we reinterpret $\forall X$ to mean for all $X$ and with answers given continuously in $X$.  Here, the key issue is how real numbers are given.  One possibility is that they are given as complete infinite objects.  In that case, computations will have to work with infinite objects directly rather than through finite codes and approximations.  There are several inequivalent approaches to that, and infinite computation is a fruitful area of research.  However, as far as it is known, humans and (physical) computers work only with finite objects, and generally, constructivists view real numbers as descriptions and as limits of finite approximations rather than as metaphysically existing infinite sets.  Through finite approximations is also how real numbers are obtained from physical systems.  Therefore, we require witnesses to provide answers based on finite segments to the sequences, which amounts to the continuity condition.  The continuity also holds in Brouwer's intuitionism and (for extensional functions) in recursive analysis.

We hope our paper makes constructive analysis more appealing and its study more productive.

\section{References}
\begin{description}
\item[{[Bishop 1967]}]  Bishop, Errett.  \textit{Foundations of Constructive Mathematics}.  McGraw-Hill, New York.

\item[{[Ishihara 2004]}]  Ishihara, Hajime.  ``Informal Constructive Reverse Mathematics.''  CDMTCS Reseearch Report Series.  Available:  \verb+http://+\\*\verb+www.cs.auckland.ac.nz/CDMTCS/researchreports/229hajime.pdf+

\item[{[Kleene 1965]}]  Stephen Kleene and Richard Vesley.  \textit{The foundations of intuitionistic mathematics, especially in relation to recursive functions}.  Amsterdam, North-Holland Pub. Co.

\item[{[Moschovakis 2005]}]  Moschovakis, Joan.  \textit{Notes on the Foundations of Constructive Mathematics}.  \\* Available:  \verb+http://www.math.ucla.edu/~joan/newnotes.ps+ 

\item[{[Simpson 1999]}] Simpson, Stephen.  \textit{Subsystems of Second Order Arithmetic}.  Springer, Berlin.

\item[{[Troelstra 1977]}]  ``Aspects of Constructive Mathematics'' in \textit{Handbook of Mathematical Logic}, edited by Jon Barwise.  
\end{description}

\end{document}